\documentclass[twocolumn]{arxiv}

\usepackage{natbib}
\usepackage{amssymb}
\usepackage{amsmath}
\usepackage{graphicx}
\graphicspath{ {fig/} }
\usepackage{epstopdf}

\usepackage{float}

\floatstyle{ruled}
\newfloat{algorithm}{thp}{lop}
\floatname{algorithm}{Algorithm}
\newfloat{subroutine}{thp}{lop}
\floatname{subroutine}{Subroutine}

\begin{document}

\begin{frontmatter}

\title{A lifting method for analyzing distributed synchronization on the unit sphere\thanksref{footnoteinfo}}

\thanks[footnoteinfo]{This work was supported by the Luxembourg National Research Fund FNR (8864515, C14/BM/8231540) and by the University of Luxembourg internal research project PPPD.}

\author[Johan]{Johan Thunberg}\ead{johan.thunberg@uni.lu}, 
\author[Johan]{Johan Markdahl}\ead{johan.markdahl@uni.lu}, 
\author[florian]{Florian Bernard}\ead{f.bernardpi@gmail.com},
\author[Johan]{Jorge Goncalves}\ead{jorge.goncalves@uni.lu}

\address[Johan]{Luxembourg Centre for Systems Biomedicine, University of Luxembourg,
    6, avenue du Swing, 
L-4367 Belvaux, Luxembourg}

\address[florian]{Max Planck Institute for Informatics,
Saarland Informatics Campus,
Campus E1 4,
66123 Saarbrucken,
Germany}


\begin{abstract}
This paper introduces a new lifting method for analyzing convergence of  continuous-time distributed synchronization/consensus systems on the unit sphere. Points on the $d$-dimensional unit sphere are lifted to the $(d+1)$-dimensional Euclidean space. The consensus protocol on the unit sphere is the classical one, where agents move toward weighted averages of their neighbors in their respective tangent planes. Only local and relative state information is used. The directed interaction graph topologies are allowed to switch as a function of time. The dynamics of the lifted variables are governed by a nonlinear consensus protocol for which the weights contain ratios of the norms of state variables. We generalize previous convergence results for hemispheres. For a large class of consensus protocols defined for switching uniformly quasi-strongly connected time-varying graphs, we show that the consensus manifold is uniformly asymptotically stable relative to closed balls contained in a hemisphere. Compared to earlier projection based approaches used in this context such as the gnomonic projection, which is defined for hemispheres only, the lifting method applies globally. With that, the hope is that this method can be useful for future investigations on global convergence.
\end{abstract}

\end{frontmatter}

\section{Introduction}\label{sec:introduction}
\noindent
This paper considers systems of agents continuously evolving on $\mathbb{S}^{d-1}$, where $d \geq 2$. The interactions between the agents are changing as a function of time. 
For such systems we are analyzing a large class of distributed synchronization/consensus control laws. The analysis tool is a lifting method, where an equivalent consensus protocol is analyzed in the ambient space that embeds the sphere. In comparison to projection methods that have been used in this context---e.g., the gnomonic projection---the proposed method is not locally but globally defined on the unit sphere. The control action is performed in the tangent plane. Only relative information between neighboring agents is used in the control laws.  Under the assumption that the time-varying graph is uniformly quasi-strongly connected, we show that the consensus manifold is globally uniformly asymptotically stable relative to any closed ball on the sphere contained in an open hemisphere. 

Synchronization on the circle, i.e., $d = 2$, is closely related to synchronization of oscillators~\citep{dorfler2014synchronization} and it is equivalent to synchronization on $\mathsf{SO}(2)$, where several applications exist such as flocking in nature and alignment of multi-robot systems.  Also for the two-dimensional sphere, i.e., $d = 3$, there are several applications such as formation flying and flocking of birds; consider for example a multi-robot system in 3D, where the relative directions between the robots are available and the goal is to align those. For higher dimensional spheres there are currently related problems such as distributed eigenvector computation, but concrete applications might arise in the future. 

The control laws at hand---and slight variations or restrictions on the graph topologies, switchings of the graphs, dimensions of the sphere, and the nonlinear weights in the control laws etc.---have been studied from various perspectives~\citep{scardovi2007synchronization,Sarlette2009,Olfati-Saber2006,Li2014,li2015collective}. There has recently been new developments~\citep{pereira2015,pereira2016,markdahl2016towards,markdahl2016}. 
In~\cite{markdahl2017almost}, almost global consensus is shown by characterization of all equilibrium points when the graph is symmetric and constant (time-invariant). It is shown that the equlibria not in the consensus manifold are unstable and the equilibra in the consensus manifold are stable. A similar technique is used in \cite{Tron2012} to show that a consensus protocol on $\mathsf{SO}(3)$ is almost globally asymptotically stable. Now, the above-mentioned results about almost global convergence come at a price. Static undirected graph topologies are assumed as well as  more restrictive classes of weights in the control protocols. Furthermore, compared to \cite{markdahl2017almost}, the right-hand sides of the system dynamics is not necessarily an intrinsic gradient and the linearization matrices at equilibriums are not necessarily symmetric. Hence, we cannot use the result due to~\cite{lojasiewicz1982trajectoires} about point convergence for gradient flows. This inspired us to take a closer look at methods that transform the consensus problem on the unit sphere (or a subset thereof) to an equivalent consensus problem in $\mathbb{R}^{d}$.
Before we address the method---referred to as a lifting method---we briefly make some connections to the related problem of consensus on $\mathsf{SO}(3)$.

The problem of consensus on $\mathsf{SO}(3)$ has been extensively studied~\citep{Sarlette2009a,Ren2010,Sarlette2010,Tron2013,Tron2014,Deng2016,johan02}. 
There is a connection between that problem and the problem of consensus on $\mathbb{S}^3$ 
when the unit quaternions are used to represent the rotations. For those, the gnomonic projection can be used to show consensus on the unit-quaternion sphere~\citep{thunberg2014distributed,thunbergaut}. In another line of research, several methods have been  introduced where control laws based on only relative information have been augmented with additional auxiliary (or estimation) variables, which  are communicated between neighboring agents. By doing so, results about almost global convergence to the consensus manifold are achieved~\citep{AS-RS:09,thunberg2017dynamic}. The latter of these two publications provides a control protocol for Stiefel manifolds, with the unit sphere and $\mathsf{SO}(d)$ as extreme cases. A similar technique had previously been used for the sphere~\citep{scardovi2007synchronization}. The idea of introducing auxiliary variables also extends to the related distributed optimization problem in~\cite{thunberg2017distributed}. In contrast to the mentioned works, in this paper we are not assuming additional communication between the agents by means of auxiliary variables. Instead only relative information is used in the protocols. In a practical setting (considering the case $d = 3$), such information can be measured by for example a vision sensor and requires no explicit communication between the agents.

In the proposed lifting method, we lift the states from the $(d-1)$-dimensional sphere into $\mathbb{R}^{d}$. The non-negative weights in the consensus protocol for the states in the lifting space are nonlinear functions. Each agent moves in a direction that is a weighted combination of the directions to the neighbors. The weights contain rational functions of the norms of the states of the agents. Since these rational functions are not well-defined at the origin, fundamental questions arise about existence, uniqueness, and invariance of sets. Those questions are answered with positive answers. 
The hope is that this lifting method will serve as a stepping-stone to future analysis on (almost) global convergence to the consensus manifold on the unit sphere. Compared to the approach in \cite{markdahl2017almost} where all the ``bad'' equlibria on $\mathbb{S}^{d-1}$ were characterized, we only need to characterize one point, which is the origin in the ``lifted space''. If we were to show that this point has a region of attraction that is of measure zero, we would have equivalently shown the desired result about almost global convergence on the unit sphere (assuming $d \geq 3$). However, the non-differentiability of this point remains  an additional challenge.

\section{Preliminaries}\label{sec:preliminaries}
\noindent
We begin this section with some \textit{set-definitions}.
The $(d-1)$-dimensional unit sphere is 
$$\mathbb{S}^{d-1} = \{y \in \mathbb{R}^{d}: \|y\|_2 = 1\}.$$
The special orthogonal group in dimension $d$ is
$$\mathsf{SO}(d) = \{Q \in \mathbb{R}^{d \times d} : Q^T = Q^{-1}, \text{det}(Q) = 1\}.$$
The set of skew symmetric matrices in dimension $d$ is 
$$\mathsf{so}(d) = \{\Omega \in \mathbb{R}^{d \times d} :  \Omega^T = - \Omega \}.$$
The set $\mathcal{H} \subset \mathbb{S}^{d-1}$ is an open \textit{hemisphere} if there is $v \in \mathbb{S}^{d-1}$ such that $\mathcal{H} = \{w \in \mathbb{S}^{d-1}: w^Tv > 0\}$. 

We consider a multi-agent system with $n$ agents.
Each agent has a corresponding state $x_i(t) \in \mathbb{S}^{d-1}$ for 
$t \in [0, \infty)$. The initial state of each agent $i$ at time $0$ is 
$x_{i0} \in \mathbb{S}^{d-1}$. Another way to represent the states of the agents is to use rotation matrices. Let $R_i(t) \in \mathsf{SO}(d)$ satisfy $R_i(t)p = x_i(t)$ for all $i$ and $t \in [0, \infty)$, where $p = [1, 0, \ldots, 0]^T$ is the \textit{north pole}; we also define $-p$ as the \textit{south pole}. Let $R_{i0}p = x_{i0}$ for all $i$, where $R_{i0}$ is the initial $R_i$-matrix at time $0$. The $R_i$-matrices can be interpreted as transformations from body coordinate frames---denoted by $\mathcal{F}_i$'s---of the agents to a world  coordinate frame $\mathcal{F}_W$. They are transforming the unit vector $p$ in the body frames to the corresponding unit vector (or point on the unit sphere) in the world coordinate frame. 
The $R_i$'s
and their dynamics are not uniquely defined, but this is not of importance for the analysis. We choose to define the dynamics of the $R_i$'s according to \eqref{eqq:3} below. 

The dynamics of the $x_i$-vectors are given by
\begin{equation}
\label{eqq:1}
\dot{x}_i = (I - x_ix_i^T)R_i[0, v_i^T]^T = R_i[0,v_i^T]^T, 
\end{equation}
where $v_i(t)\in \mathbb{R}^{d-1}$ for all $t$. The $v_i$-vectors are the controllers for the agents and those are defined in the body coordinate frames, i.e., the $\mathcal{F}_i$'s. For the $R_i$-matrices the dynamics is 
\begin{align}
\label{eqq:3}
\dot{R}_i & = R_i\begin{bmatrix}
0 & -v_i^T \\
v_i & 0
\end{bmatrix}.
\end{align}
The matrix on the right-hand side of $R_i$
in \eqref{eqq:3} is an element of $\mathsf{so}(d)$. 
The control is performed
in the tangent space of the sphere, which means that there are $d-1$ 
degrees of freedom for the control. This is the reason why the  $v_i$-vectors are $(d-1)$-dimensional. Before we proceed, we provide some additional explanation for the expression in the right-hand side of \eqref{eqq:3}. According to its definition, the first column of $R_i$ is equal to $x_i$ and by multiplying $\dot x_i$ by $R_i^T$ from the right we obtain---due to \eqref{eqq:1}---the following expression
$$R_i^T\dot x_i = [0, v_i^T]^T.$$
This means that $$R_i^T\dot R_i = \begin{bmatrix}
0 & \star \\
v_i & \star
\end{bmatrix},$$
where the $\star$-parts are left to be chosen. We know that the matrix in the right-hand side above needs to be skew symmetric, since $R_i$ is a rotation matrix. We also know that the first column of it must be equal to $[0, v_i^T]^T$. The matrix of minimum Euclidean norm that fulfills these two requirements is equal to $$\begin{bmatrix}
0 & -v_i^T \\
v_i & 0
\end{bmatrix},$$
i.e., the one we chose in the right-hand side of \eqref{eqq:3}.

We will study a class of distributed synchronization/consensus control laws on the unit sphere, where the agents are moving in 
directions comprising conical combinations of directions to neighbors. In this protocol only local and relative information is used. Before we 
provide these control laws we introduce directed graphs and time-varying directed graphs. 

A directed graph $\mathcal{G}$ is a pair $(\mathcal{V}, \mathcal{E})$, where $\mathcal{V} = \{1, \ldots, n\}$ is the node-set and $\mathcal{E} \subset \mathcal{V} \times \mathcal{V}$ is the edge-set. Each node in the node set corresponds to a unique agent.  
The set $\mathcal{N}_i \in
\mathcal{V}$ is the neighbor set or neighborhood 
of agent $i$, where $j \in \mathcal{N}_i$
if and only if $(i,j) \in \mathcal{E}$. We continue with the
following definitions addressing connectivity of directed graphs.

In a directed graph $\mathcal{G}$, a \textit{directed path} is a sequence of distinct nodes, such that any consecutive pair of nodes in the sequence comprises an edge in the graph. 
We say that $i$ is \textit{connected to} $j$ if there is a directed path from $i$ to $j$. 
We say that the graph is \textit{quasi-strongly connected} if there is at least one node that is a center or a root node in the sense that all the other nodes are connected to it.  We say that the graph is \textit{strongly connected} if for all $(i,j) \in \mathcal{V} \times \mathcal{V}$ it holds that $i$ is connected to $j$. 

Now we define \textit{time-varying graphs}. We define those by first defining \textit{time-varying neighborhoods}. The time-varying neighborhood $\mathcal{N}_i(t)$ of agent $i$
is a piece-wise constant right-continuous set-valued function that maps from $\mathbb{R}$ to $2^{\mathcal{V}}$. We assume that there is $\tau_D >0$ such that 
$\inf_{k}(\tau_{i({k+1})} - \tau_{ik}) > \tau_D$ for all $i$, where  $\{\tau_{ik}\}_{k = -\infty}^{\infty}$ is the set of time points of discontinuity of $\mathcal{N}_i(t)$. The constant $\tau_D$ is as a lower bound on the dwell-time between any two consecutive switches of the topology.
We define the \textit{time-varying graph} $\mathcal{G}(t) = (\mathcal{V},\mathcal{E}(t))$ as
 $$\mathcal{G}(t) = (\mathcal{V},\mathcal{E}(t)) = (\mathcal{V},\bigcup_i\bigcup_{j \in \mathcal{N}_i(t)}\{(i,j)\}).$$
Furthermore, the \textit{union graph} of $\mathcal{G}(t)$ during the time
interval $[t_1,t_2)$ is defined by
\begin{equation*}
\mathcal{G}([t_1, t_2))
= \textstyle\bigcup_{t\in[t_1, t_2)} \mathcal{G}(t)
= (\mathcal{V},\textstyle\bigcup\nolimits_{t\in[t_1, t_2)}\mathcal{E}(t)),
\end{equation*}
where $t_1 < t_2 \leq \infty$.
We say that the graph $\mathcal{G}(t)$ is \textit{uniformly
(quasi-) strongly connected} if there exists a constant $T >0$ such that the
union graph $\mathcal{G}([t, t + T))$ is (quasi-) strongly connected for all $t$.

Now we provide the synchronization protocol to be studied. For each agent $i$, the controller is $v_i$ is defined by
\begin{align}
\label{eqq:4}
\begin{bmatrix}
0 \\
v_i
\end{bmatrix} & = \begin{bmatrix}
0 & 0 \\
0 & I_{d-1}
\end{bmatrix} \sum_{j \in \mathcal{N}_i(t)}f_{ij}(\|x_{ij} - p\|)x_{ij}, 
\end{align}
where $x_{ij} = R_i^Tx_j =  R_i^TR_jp$, which is $x_j$ represented in the frame $\mathcal{F}_i$. The $x_{ij}$'s are what we refer to  as relative information and the control law \eqref{eqq:4} is constructed by only such information. 
For each $(i,j)$, it holds that $f_{ij}: \mathbb{R} \rightarrow \mathbb{R}$. 
The $f_{ij}$-functions are assumed to be Lipschitz and attain positive values for positive arguments. The $\mathcal{N}_i(t)$'s are neighborhoods of a time-varying directed graph $\mathcal{G}(t)$, whose connectivity is at least uniformly quasi-strong. These control laws will be analyzed in the paper. 

The expressions in \eqref{eqq:4} are more easily understood if they are expressed in the world frame $\mathcal{F}_{W}$. We define
\begin{align}
\label{eqq:5}
u_i = (I - x_ix_i^T) \sum_{j \in \mathcal{N}_i(t)}f_{ij}(\|x_j - x_i\|)(x_j - x_i), 
\end{align}
for all $i$, 
which is $[0, v_i^T]^T$ expressed in the frame $\mathcal{F}_W$. The vector $u_i$ is the sum of the positively weighted directions to the neighbors of agent $i$, projected onto the tangent space at the point $x_i$. Also for analysis purposes, \eqref{eqq:5} is easier to work with than \eqref{eqq:4}. The closed loop system is
\begin{align}
\label{eqq:6}
\dot{x}_i & = (I - x_ix_i^T)\sum_{j \in \mathcal{N}_i(t)}f_{ij}(\|x_j - x_i\|)(x_j - x_i), \\
\nonumber
 & = (I - x_ix_i^T)\sum_{j \in \mathcal{N}_i(t)}f_{ij}(\|x_j - x_i\|)x_j,
\end{align}
for all $i$.  

Let $x = [x_1^T, x_2^T, \ldots, x_n^T]^T$ and $x_0 = [x_{10}^T, x_{20}^T, \ldots, x_{n0}^T]^T$.  
We define the set 
$$\mathcal{A} = \{x: x_i = x_j \text{ for  all }i,j\},$$
which is the synchronization/consensus set. Throughout the paper we assume that the closed-loop dynamics of the 
system is given by \eqref{eqq:6}.
We study the convergence of $x(t)$ to the consensus set $\mathcal{A}$. When we talk about convergence we refer to the concepts below. 

For the system \eqref{eqq:6}, we say that the set $\mathcal{A}$ is \textit{attractive} relative to a forward invariant set $\mathcal{S} \subset (\mathbb{S}^{d-1})^n$ if $$(x_{0} \in \mathcal{S}) \Longrightarrow (\text{dist}(x(t), \mathcal{A}) \rightarrow 0\text{ as } t \rightarrow \infty)$$
where $\text{dist}(v,\mathcal{A}) = \inf_{w \in \mathcal{A}}\|w - v\|$. Furthermore, we say that the set $\mathcal{A}$ is \textit{globally uniformly asymptotically stable} relative to a forward invariant compact set $\mathcal{S} \subset (\mathbb{S}^{d-1})^n$ if 
\begin{enumerate}
\item for every $\epsilon > 0$ there is $T(\epsilon) > 0$ such that 
$(x_0 \in \mathcal{S}) \Longrightarrow (\text{dist}(x(T(\epsilon)),\mathcal{A}) \leq \epsilon);$ \\
\item for every $\epsilon >0$ there is $\delta(\epsilon) > 0$ such that \\
$(x_0 \in \mathcal{S} \text{ and } \text{dist}(x_0,\mathcal{A}) \leq \delta) \Longrightarrow (\text{dist}(x(t),\mathcal{A}) \leq \epsilon$
for all $t \geq 0$). 
\end{enumerate}
The equivalent definitions to the above will also be used (after changing the sets $\mathcal{S}$ and $\mathbb{S}^{d-1}$) for other systems evolving in $(\mathbb{R}^d)^n$
or linear subspaces thereof. \textit{Forward invariance}, or simply \textit{invariance}, of a set means that if the initial state is contained in the set, then the state is contained in the set for all future times.  	

The two concepts of global convergence respective almost global convergence relative to a forward invariant set $\mathcal{S}$ refer to, respectively, the situations where convergence occur for all initial points in $\mathcal{S}$ and convergence occur for all initial points in a set $\mathcal{B}$ where $\mathcal{S} - \mathcal{B}$ has measure zero.

\section{Projection methods}
\noindent
Before we continue to present the lifting method, we show how projection 
based methods can be used to analyze consensus on hemispheres. In particular we consider two such methods. 
The two methods are such that the $x_i$-vectors are projected down onto a $(d-1)$-dimensional linear subspace of $\mathbb{R}^{d}$. The symbol $y_i$ is used to denote the projection variable for $x_i$ in both methods.

\subsection{Equatorial plane projection}
\noindent
The equatorial plane projection simply projects all the states onto a $(d-1)$-dimensional hyperplane (that contains the origin). This plane separates the sphere into two hemispheres. If all the agents are positioned on one of those hemispheres, one can easily show that they reach consensus provided that the graph has strong connectivity. This projection method is appealing because the projections are simple and the convergence proof is straightforward. It is interesting that results from the literature about convergence on hemispheres (and slightly more general ones where the graph is assumed to be time-varying) can easily be shown with this simple projection. 

Now, formally, the $x_i$-states are projected onto the equatorial plane whose normal is equal to $p = [1,0, \ldots, 0]^T$ in the world coordinate frame $\mathcal{F}_W$. 

\begin{figure}[ht]
  \centering
  \includegraphics[width=\columnwidth]{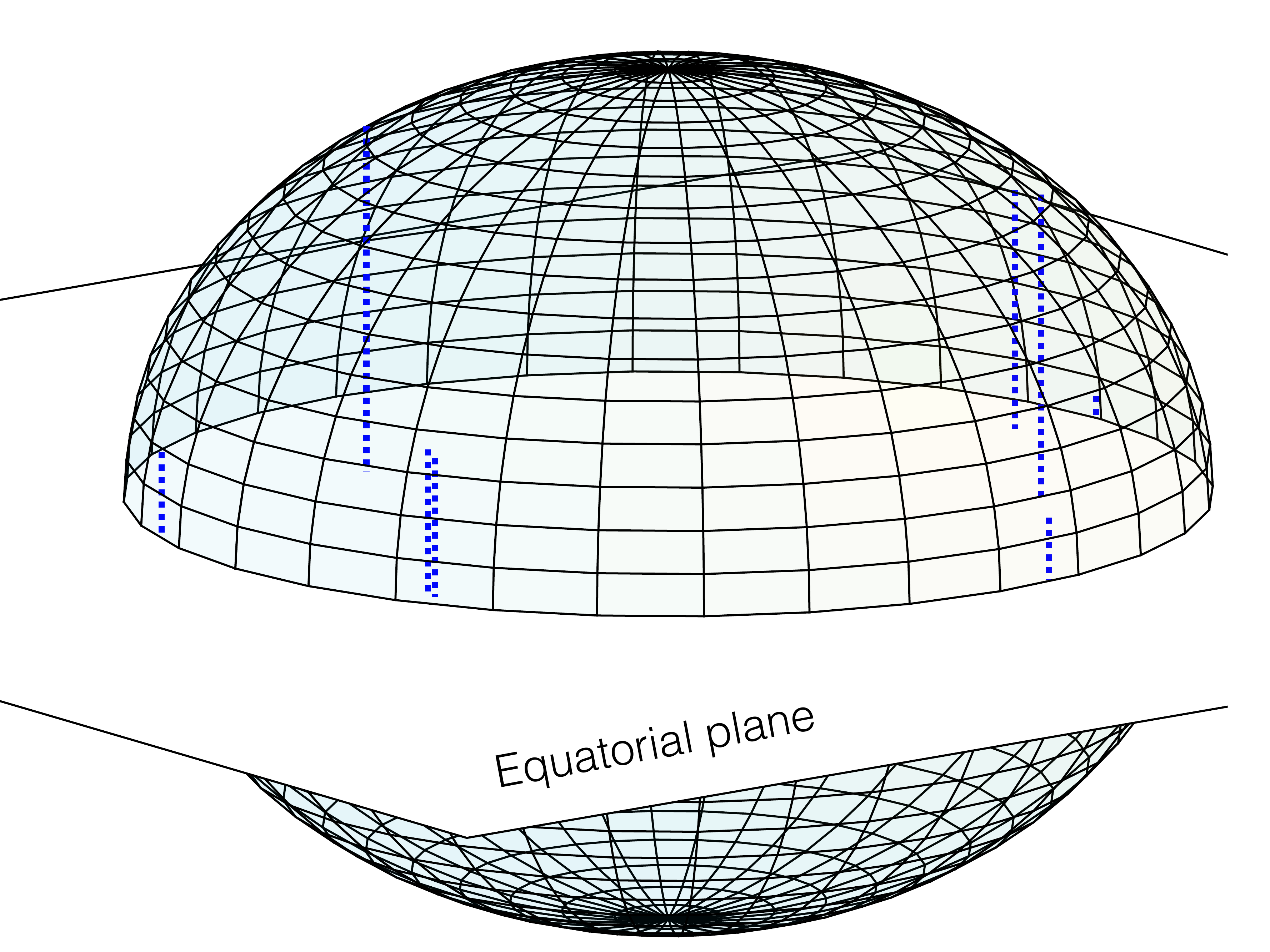}
  \caption{Illustration of the equatorial plane projection.}
  \label{fig:1}
\end{figure}

The projected state $y_i$ is defined by 
\begin{equation}
\begin{bmatrix}
0 \\
y_i
\end{bmatrix} = P_{\text{equ}}x_i = \begin{bmatrix}
0 & 0 \\
0 & I_{d-1}
\end{bmatrix}x_i.
\end{equation}
This is illustrated in Fig.~\ref{fig:1} for the dimension $d = 3$. Points on the northern hemisphere, i.e., the $x_i$'s satisfying $p^Tx_i > 0$, are projected down onto the equatorial plane. For each point there is a blue dotted line between the point and its projection.

On the northern hemisphere $x_i \mapsto y_i$ is a diffeomorphism. The mapping $y_i \mapsto x_i$
is defined by
\begin{align}
[x_i]_1 & = \sqrt{1 - \|y_i\|^2}\\
[x_i]_k & = [y_i]_k, \text{ for } k\geq 2,
\end{align}
where $[x_i]_k$ and $[y_i]_k$ are the $k$'th elements of $x_i$ and $y_i$, respectively.

By using this projection one obtains a local  convergence result for hemispheres. 

\begin{prop}\label{prop:1}
Suppose controller \eqref{eqq:4} is used for each agent $i$ and suppose $\mathcal{G}(t)$ is uniformly strongly connected. If $p^Tx_{i0} > 0$ for all $i$, it holds that $\mathcal{A}$ is attractive for the closed loop system \eqref{eqq:6}. 
\end{prop}

\emph{Proof:} \quad
We will use Theorem 1 in \cite{Johan_lyap_2017}. Under the condition that the graph is uniformly strongly connected, if we can show that any closed disc (or ball) in the equatorial place with radius less than $1$ is forward invariant for the $y_i$'s and we can find a function $V: \mathbb{R}^{d-1} \rightarrow \mathbb{R}^+$ such that $V$ is 1) positive definite, 2) $\max_{i \in \mathcal{V}}V(y_i(t))$ is decreasing as a function of $t$, and 3) $\dot V(y_i(t))$ is strictly negative if $i \in \arg\max_{j \in \mathcal{V}}\{V(y_j)\}$ and there is $j \in \mathcal{N}_i(t)$ such that $y_i \neq y_j$. Then the $y_i$'s converge to a consensus formation. This in turn implies, since $y_i \mapsto x_i$ is a diffeomorphism, that the $x_i$'s converge to a consensus formation, i.e., the set $\mathcal{A}$ is attractive. 

Let $V = \|\cdot\|^2$, i.e., $V(y_i) = y_i^Ty_i$, which obviously satisfy condition 1). For $\|y_i\| < 1$ it holds that
\begin{align}
\nonumber
& \quad \: \dot{V}(y_i) \\
\nonumber
 & = \sum_{j \in \mathcal{N}_i(t)} g_{ij}(y_i, y_j)\bigg (\left(\sqrt{1 - \|y_i\|^2}\sqrt{1 - \|y_i\|^2}\right )y_i^Ty_j \\
\label{eqq:8}
& \quad \:- \left(\sqrt{1- \|y_i\|^2} \sqrt{1- \|y_{\smash{j}}\|^2}\right )\|y_i\|^2 \bigg),
\end{align} 
where $g_{ij}(y_i, y_j) = \\ f_{ij}(\|[\sqrt{1- \|y_j\|^2} - \sqrt{1- \|y_i\|^2}, (y_j - y_i)^T]^T\|)$. 

Now, at time $t$, assume $i \in \text{arg}\max_{j \in \mathcal{V}}\{V(y_j)\}$ and assume $\|y_i\| < 1$. The following observations imply that conditions 2) and 3) hold. It holds that $g_{ij}(y_i,y_j) \geq 0 $ and the inequality is strict if $y_j \neq y_i$. It holds that $\|y_i\|^2 \geq y_j^Ty_j$ and the inequality is strict if $y_j \neq y_i$. It holds that 
$$\left(\sqrt{1- \|y_i\|^2} \sqrt{1- \|y_{\smash{j}}\|^2}\right ) \geq \left(\sqrt{1 - \|y_i\|^2}\right )^2 \geq 0.$$

We also see that any closed disc with radius less than $1$ 
is forward invariant for the $y_i$'s. 
\hfill $\blacksquare$

By a change of coordinates, we obtain the following generalization. 
\begin{cor}\label{cor:1}
Suppose controller \eqref{eqq:4} is used for each agent $i$ and suppose $\mathcal{G}(t)$ is uniformly strongly connected. If all the $x_i$'s are contained in an open hemisphere it holds that $\mathcal{A}$ is attractive for the closed loop system \eqref{eqq:6}. 
\end{cor}

A main problem, with the equatorial plane projection is that the convex hull of the projected variables is not necessarily forward invariant. This means that the projected variables are not following a consensus protocol. This is also the reason why we settle for the strong connectivity assumption about the graph, i.e., that it is uniformly strongly connected. However, the projected variables under the gnomonic projection---introduced in the subsequent section---do follow a consensus protocol, which, in turn, allows for more general convergence results.

\subsection{The gnomonic projection}
\begin{figure}[ht]
  \centering
   \includegraphics[width=\columnwidth]{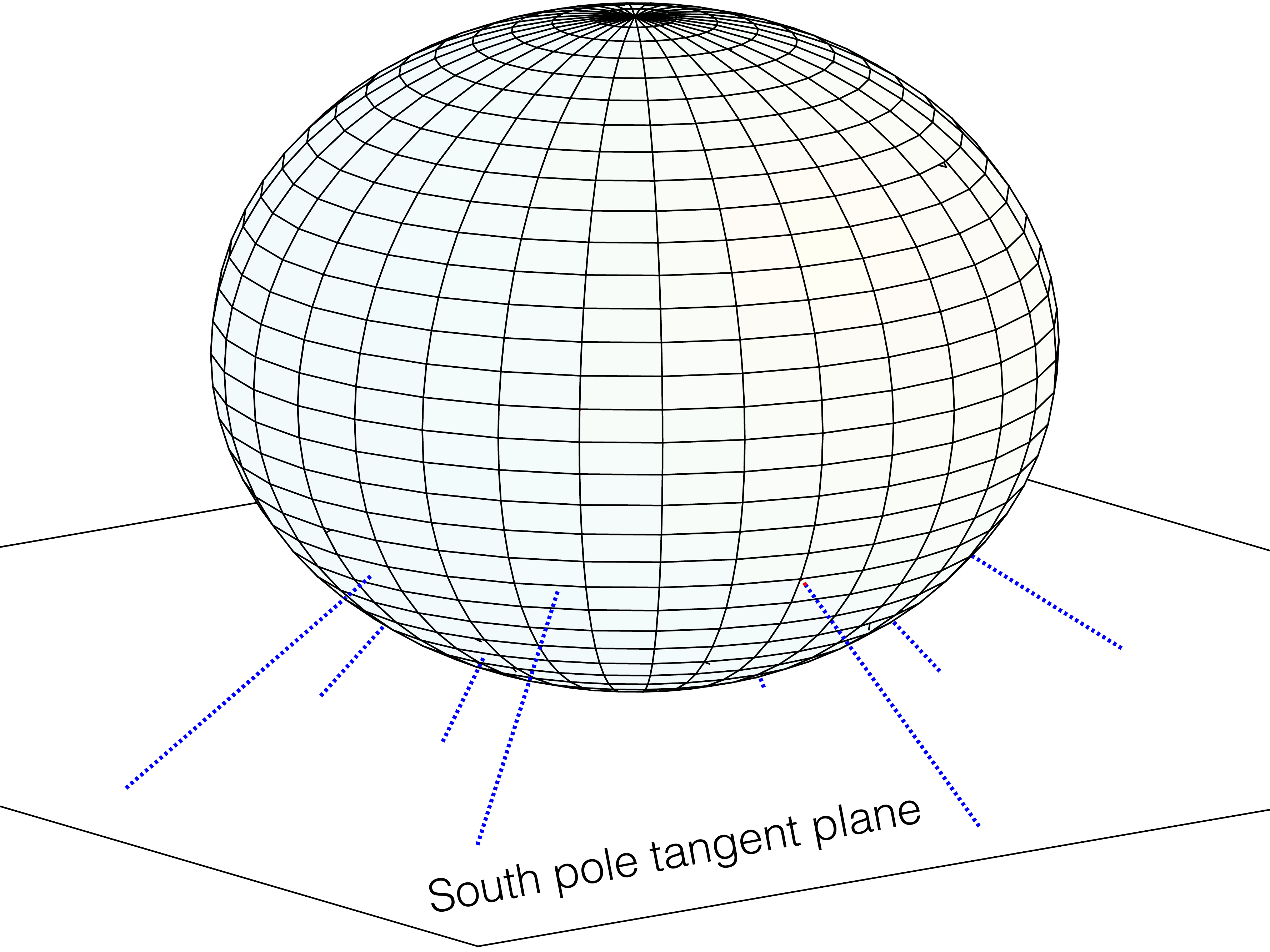}
  \caption{Illustration of the gnomonic projection.}
  \label{fig:2}
\end{figure}

\noindent
The gnomonic projection projects an open hemisphere onto a tangent plane at a point on the sphere. We will use the convention of projecting the points on the southern hemisphere defined as $\{x \in \mathbb{S}^{d-1}: p^Tx < 0\}$ onto the tangent plane at the south-pole, i.e., at the point $-p$. The projection of $x_i$ is the intersection between the tangent plane and the line that passes through the origin and $x_i$. This projection is illustrated in Fig.~\ref{fig:2}, where several points are projected. 
The gnomonic projection has the property that segments  of great circles on the sphere (geodesics) correspond to straight line segments in the projection plane. 
One can show that a consensus algorithm on the open hemisphere corresponds to a consensus protocol for the projected states. It should be emphasized that the gnonomic projection method is not new. It is claimed to have been invented by the Greek philosopher Thales of Miletus somewhere around 624--546 BCE~\citep{alsina2015mathematical}. Its first appearance in a subject related to the one addressed in this paper, was probably in \cite{Hartley2011} and subsequently in  \cite{hartley2013rotation} in the context of rotation averaging.
Later the gnomonic projection was used as a tool to show consensus on the open hemisphere~\citep{thunberg2014distributed,johan02}. In those latter works, the three-dimensional (unit-quaternion) sphere was considered in the context of attitude synchronization. Recently, the gnomonic projection was also  considered for arbitrary dimensions~\citep{lageman2016consensus}. It should be emphasized that the graph was not time-varying in that context.

Formally we define $y_i$,  the projection of $x_i$, by 
\begin{equation}
\begin{bmatrix}
-1 \\
y_i
\end{bmatrix} = \frac{1}{|[x_i]_1|}x_i, 
\end{equation}
which is a diffeomorphism from the open southern hemisphere to the tangent plane of the south pole.
Suppose controller
\eqref{eqq:4} is used, i.e., the closed loop dynamics is given by \eqref{eqq:6}. If $p^Tx_i < 0$ for all $i$, i.e., all the $x_i$'s are located on the southern hemisphere, it holds that the dynamics of the $\dot{y}_i$'s is on the form
\begin{align}
\label{eqq:nisse:1}
\dot{y}_i = \sum_{j \in \mathcal{N}_i(t)}h_{ij}(y)(y_j - y_i),
\end{align}
where $y = [y_1^T, y_2^T, \ldots, y_n^T]^T$, and it can be shown that the $h_{ij}(y)$'s are locally Lipschitz and globally Lipschitz on any bounded set. This can be used (as an alternative to using the lifting method) to prove the result in Proposition~\ref{prop:10} in next Section, which is stronger than that in Corollary~\ref{cor:1}.

\section{The lifting method}
\begin{figure}[ht]
  \centering
  \includegraphics[width=\columnwidth]{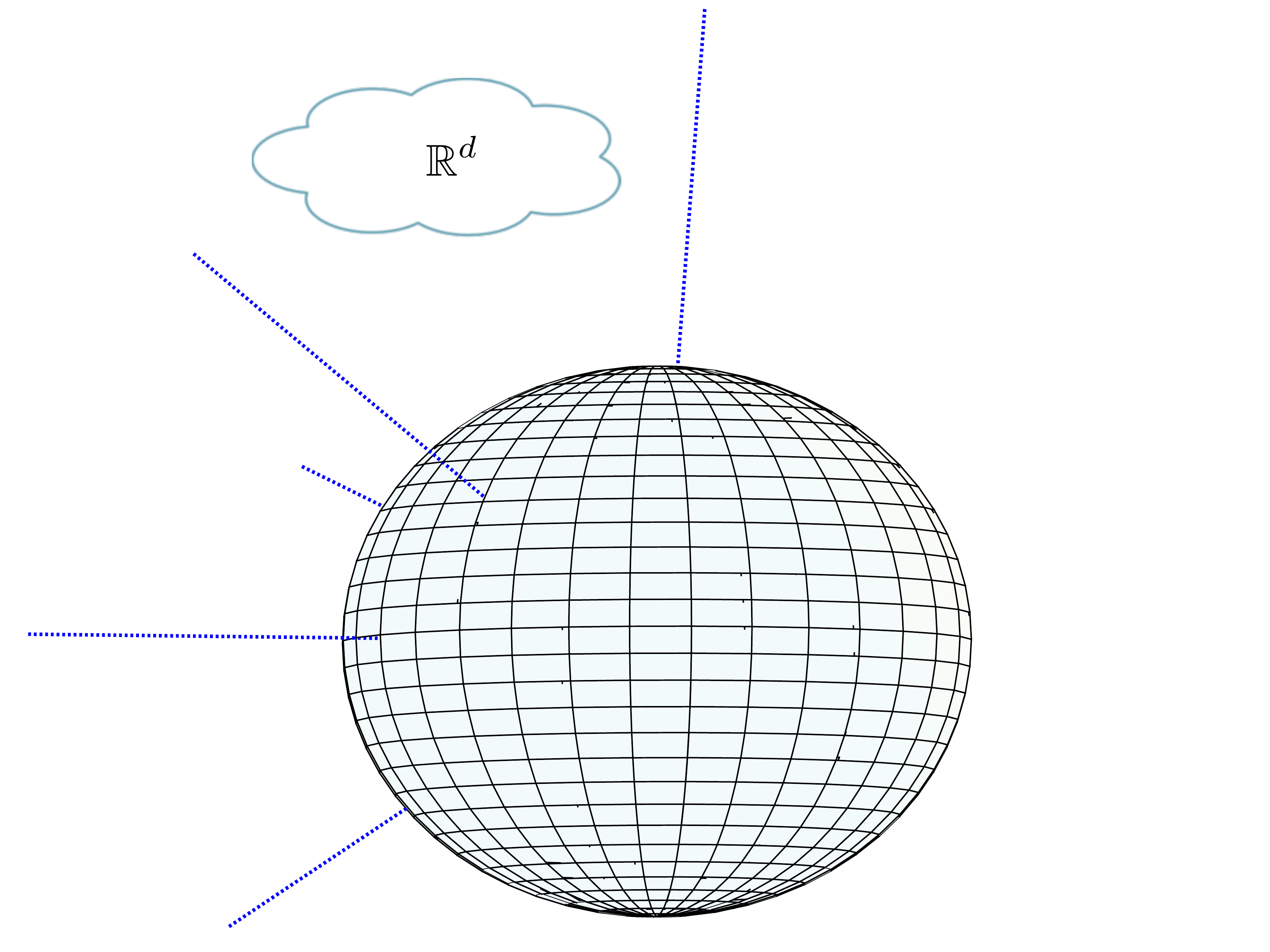}
  \caption{Illustration of the lifting method.}
  \label{fig:3}
\end{figure}

\noindent
In this section we propose a method where the $x_i$'s are not projected onto a ($d-1$)-dimensional plane, but rather relaxed to be elements in $\mathbb{R}^{d}$. Those elements, we call them $z_i$'s, are then projected down onto the sphere $\mathbb{S}^{d-1}$ to create the $y_i$'s (which in this case are equivalent to the $x_i$'s). The method can thus be seen as the inverse procedure to the two in the previous section. Provided $z_i \neq 0$, the projection is given by 
\begin{equation}
y_i = \frac{z_i}{\|z_i\|}.
\end{equation}
This projection as well as the lifting is illustrated in Figure~\ref{fig:3}. Points in $\mathbb{R}^d$ are projected down onto the sphere in the sense of minimizing the least squares distance.

We let $z_i(t) \in \mathbb{R}^{d}$ be governed by the following dynamical system
\begin{align}
\label{eqq:11}
\dot{z}_i = \sum_{j \in \mathcal{N}_i(t)}f_{ij}\left(\left\|\frac{z_j}{\|z_j\|} - \frac{z_i}{\|z_i\|}\right\|\right)\frac{\|z_i\|}{\|z_{\smash{j}}\|}(z_j - z_i), 	
\end{align}
for all $i$.
 Let the initial state of the system be $z_0 = [z_{10}^T, z_{20}^T, \ldots, z_{n0}^T]^T$. 
Equation \eqref{eqq:11} describes a consensus protocol with nonlinear weights that contain rational functions of the norms of the states.
The question is how this dynamical system is related to \eqref{eqq:6}. The following proposition provides the answer. 

\begin{prop}\label{prop:5}
 Suppose that all the $z_{i0}$'s are not equal to zero. On the time interval $[0, \infty)$ the dynamics for the $y_i$'s is given by 
\begin{equation}\label{eqq:301}
\dot{y}_i = (I - y_iy_i^T)\sum_{j \in \mathcal{N}_i(t)}f_{ij}(\|y_j - y_i\|)(y_j - y_i),
\end{equation}
i.e., it is the same as \eqref{eqq:6}.
\end{prop}

\emph{Proof:} \quad
Proposition~\ref{prop:3} below provides the result that the solution to \eqref{eqq:11} is well-defined on the interval $[0, \infty)$ and $(\mathbb{R}^d\backslash \{0\})^n$ is forward invariant. Given that result, the $y_i$'s and their derivatives are well defined. Now, 
\begin{align*}
\dot{y}_i & = \frac{1}{\|z_i\|}(I - y_iy_i^T)\dot{z}_i \\
& = (I - y_iy_i^T)\sum_{j \in \mathcal{N}_i(t)}f_{ij}(\left\|y_j - y_i\right\|)\left(y_j - y_i\frac{\|z_i\|}{\|z_j\|}\right) \\
& = (I - y_iy_i^T)\sum_{j \in \mathcal{N}_i(t)}f_{ij}(\|y_j - y_i\|)(y_j - y_i). \hfill \quad \quad \blacksquare
\end{align*}

\begin{prop}\label{prop:3}
Suppose the dynamics for the $z_i$'s is governed by \eqref{eqq:11}. Suppose there is no $i$ such that $z_{i0} = 0$. Let $H(t,z_0)$ be the convex hull of the $z_{i}(t)$'s during $[0, t_f)$ when the initial condition is $z_0$. Let $t_f$ be such that the solution exists during $[0,t_f)$.
Then the solution to \eqref{eqq:11} exists and is unique for all times $t > 0$, $(H(t,z_0))^n$ is forward invariant, and the set  $(\mathbb{R}^d\backslash \{0\})^n$ is forward invariant.
\end{prop}

\emph{Proof:} \quad
We first address the claim that $(H(t,z_0))^n$ is forward invariant. It suffices to verify that for each $i$, the right-hand side of \eqref{eqq:11} is either inward-pointing relative to $H(t,z_0)$, or equal to $0$. Now, due to the structure of \eqref{eqq:11}, this is true. 

Now we address the invariance of $(\mathbb{R}^d\backslash \{0\})^n$, and, by doing that, obtain the existence and uniqueness result for the solution during $[0, \infty)$ for free, since the right-hand side of \eqref{eqq:11} is locally Lipschitz on $(\mathbb{R}^d\backslash \{0\})^n$. Now, suppose there is $i_1 \in \mathcal{V}$ and a finite time $t_1 > 0$ such that $\lim_{t \uparrow t_1} z_{i_1}(t) = 0$ and there is no $j$ and $t_0 \in [0,t_1)$ such that $\lim_{t \uparrow t_0} z_{j}(t) = 0$. This means that there is a first finite time $t_1$ for which at least one state, $z_{i_1}$ that is, attains the value $0$. The  assumption is equivalent to assuming that  $(\mathbb{R}^d\backslash \{0\})^n$ is not forward invariant. 

For $i \in \mathcal{V}$ and $t \in [0, t_1)$ it holds that 
\begin{align*}
\frac{d}{dt}{\|z_i\|} & = \sum_{j \in \mathcal{N}_i(t)}g(z_i, z_j)(\|z_i\|\theta_{ij} - \frac{~\|z_i\|^2}{\|z_j\|}), \text{ for all } i,
\end{align*}
where $g(z_i,z_j) = f(\|\frac{z_j}{\|z_j\|} - \frac{z_i}{\|z_i\|}\|)$
and $\theta_{ij} = \frac{z_i^Tz_j}{\|z_i\|\|z_j\|}$. For $z_i, z_j \neq 0$, it also holds that 
\begin{align*}
 \frac{d}{dt}\frac{\|z_i\|}{\|z_j\|} & = \sum_{k \in \mathcal{N}_i(t)}g_{ik}(z_i,z_k)\left(\frac{\|z_i\|}{\|z_j\|}\theta_{ik} - \frac{~\|z_i\|^2}{\|z_j\|\|z_k\|}\right) \\
 &\quad - \sum_{l \in \mathcal{N}_j(t)}g_{jl}(z_j,z_l)\left(\frac{\|z_i\|}{\|z_j\|}\theta_{jl} - \frac{\|z_i\|}{\|z_l\|}\right).
 \end{align*} 
We define $v_{kl} = \frac{\|z_k\|}{\|z_l\|}$ for all $k,l$ and write the 
equation above as 
 \begin{align*}
 \dot{v}_{ij} & = \sum_{k \in \mathcal{N}_i(t)}g_{ik}(z_i,z_k)(v_{ij}\theta_{ik} - v_{ij}v_{ik}) \\
 & \quad - \sum_{l \in \mathcal{N}_j(t)}g_{jl}(z_j,z_l)(v_{ij}\theta_{jl} - v_{il}).
\end{align*}
Let $\alpha > 0$ be an upper bound for the $f_{ij}$'s,
which is equivalent to an upper bound for the $g_{ij}$'s. Such a bound must exist (since the set $\mathbb{S}^{d-1} \times \mathbb{S}^{d-1}$ is compact and the function $(f_{ij} \circ \text{dist})$ is continuous on $\mathbb{S}^{d-1} \times \mathbb{S}^{d-1}$, where $\text{dist}(\cdot,\cdot)$ is the function that returns the Euclidean distance between two points in $\mathbb{R}^d$). 

Let $V(t) = \max\limits_{(i,j) \in \mathcal{V} \times \mathcal{V}}v_{ij}(t)$. On $[0, t_1)$ it holds that 
\begin{align}
\label{eqq:12}
D^+{V} \leq  3\alpha n V,
\end{align} 
where $D^+$ is the upper Dini-derivative. By using the Comparison Lemma for \eqref{eqq:12}, we can conclude that $V$ is bounded from above by $e^{3\alpha nt_1}V(0)$ on $[0,t_1)$

Now, for $i_1$ and $t \in [0, t_1)$ it holds that 
\begin{align*}
\frac{d}{dt}{\|z_{ i_1}\|}  & = \sum_{j \in \mathcal{N}_{i_1}(t)}g(z_{i_1}, z_j)(\|z_{ i_1}\|\theta_{ij} - \frac{~\|z_{ i_1}\|^2}{\|z_j\|}) \\
& \geq -n\alpha(e^{3\alpha nt_1}V(0) + 1)\|z_{i_1}\|.
\end{align*}
By using the Comparison Lemma, we can conclude that $\|z_{ i_1}(t)\| \geq \|z_{i_1 0}\|e^{-n\alpha(e^{3\alpha nt_1}V(0) + 1)}$. But this, in turn, means that $\lim_{t \uparrow t_1}\|z_{i_1}(t)\| > 0$, which is a contradiction.
\hfill $\blacksquare$

In the following proposition we make use of $H(t,z_0)$, which was defined in Proposition~\ref{prop:3}.

\begin{prop}\label{prop:7}
Suppose the dynamics for the $z_i$'s are governed by \eqref{eqq:11} and $\mathcal{G}(t)$ is uniformly quasi-strongly connected. Suppose $0 \in \mathbb{R}^d$ is not contained in the convex hull $H(t,w)$ for $w \in \mathbb{R}^{nd}$. Then the consensus set 
$\mathcal{A}_z$---defined as the set where all the $z_i$'s are equal in $(H(0,w))^n$---is globally uniformly asymptotically stable relative to $(H(0,w))^n$. Furthermore, there is a point $\bar z \in \mathbb{R}^d$ that all the $z_i$'s converge to.
\end{prop}

\emph{Proof:} \quad
Invariance of $(H(t,z_0))^n$ is an indirect consequence of the 
fact that \eqref{eqq:11} is a consensus protocol. On this set the right-hand side of \eqref{eqq:11} is Lipschitz continuous in $z$ and piece-wise continuous in $t$. Now the procedure in the rest of the proof is 
analogous to the one in Proposition~\ref{prop:10}. Since the right-hand side of \eqref{eqq:11} is Lipschitz continuous in $z$ and piece-wise continuous in $t$, we can use Theorem 2 in \cite{Johan_lyap_2017} to find a continuously differentiable function ${W}: \mathbb{R}^{d} \times \mathbb{R}^{d} \rightarrow \mathbb{R}^+$ such that 1) $\max_{(i,j) \in \mathcal{V} \times \mathcal{V}}W(z_i(t),z_j(t))$ is decreasing as a function of $t$; and 2) $\dot{W}(z_i(t),z_j(t))$ is strictly negative if $(i,j) \in \text{arg}\max_{(i,j) \in \mathcal{V} \times \mathcal{V}}W(z_i(t),z_j(t))$ and there is $k \in \mathcal{N}_i(t)$ such that $y_i \neq y_k$ or there is $l \in \mathcal{N}_j(t)$ such that $z_j \neq z_l$. The existence of such a function guarantees that the consensus set $\mathcal{A}_z$ is globally uniformly asymptotically stable relative to $(H(0,z_0))^n$. It holds that the function $\|z_i - z_j\|^2$ is such a $W$-function. Convergence to a point for all the $z_i$'s can be shown by using the facts that $(H(t,z_0))^n$ is forward invariant for all $t$ and $z$ converges to $\mathcal{A}_z$. 
\hfill $\blacksquare$

As a remark to the previous proposition, we should add that more restrictive results about attractivity of $\mathcal{A}$ can be shown by using the results in~\cite{shi2009global,lin2007state}.

\begin{prop}\label{prop:10}
Suppose the graph $\mathcal{G}(t)$ is uniformly quasi-strongly connected. Then for any closed ball $B$ contained in the hemisphere, the consensus set $\mathcal{A}$ is globally uniformly asymptotically stable relative to $B^n$ under \eqref{eqq:6}.
\end{prop}

\emph{Proof:} \quad
Forward invariance holds for $B^n$ due to the structure of the right-hand side of \eqref{eqq:6}.
Let $z_0 = x_0$.
Due to Proposition \ref{prop:7} we know that the consensus set 
$\mathcal{A}_z$ is globally uniformly asymptotically stable relative to $(H(0,z_0))^n$ and that there is a point $\bar z \in \mathbb{R}^d$ that all the $z_i$'s converge to. We also know that the projected $y_i$-variables follow the protocol \eqref{eqq:301}, which is the same as \eqref{eqq:6}. The norms of the $z_i$'s are uniformly bounded on $(H(0,z_0))^n$ and $(H(0,z(T)))^n$ is forward invariant for all $T > 0$. Thus the desired result readily follows. 
\hfill $\blacksquare$

\begin{cor}
Suppose the dynamics for the $z_i$'s is governed by \eqref{eqq:11} and suppose the graph $\mathcal{G}(t)$ is quasi-strongly connected. If the $z_i$'s converge to a point $\bar{z} \in \mathbb{R}^d$ that is not equal to zero, then the $y_i$'s converge to a point $\bar{y} \in \mathbb{S}^{d-1}$. Furthermore, if the convex hull of the $z_{i0}$'s does not contain the point zero, 
then the $y_i$'s converge to a point $\bar{y} \in \mathbb{S}^{d-1}$.
\end{cor}

\emph{Proof:} \quad
Straightforward application of Proposition~\ref{prop:5}, Proposition~\ref{prop:3}, and Proposition~\ref{prop:10}.
\hfill $\blacksquare$

Variations of Proposition~\ref{prop:10} has appeared in the literature before. The idea of using the gnomonic projection to show consensus on the hemisphere was used in~\cite{thunberg2014distributed,johan02} where restricted versions of Proposition~\ref{prop:10} were given for the dimension $d = 4$ in the context of attitude synchronization. Recently the attractivity of $\mathcal{A}$ relative to open hemispheres was established under quasi-strong graph connectivity~\citep{lageman2016consensus} using the gnomonic projection. The graph was not time-varying in that context.

\begin{figure}[ht]
  \centering
  \includegraphics[width=\columnwidth]{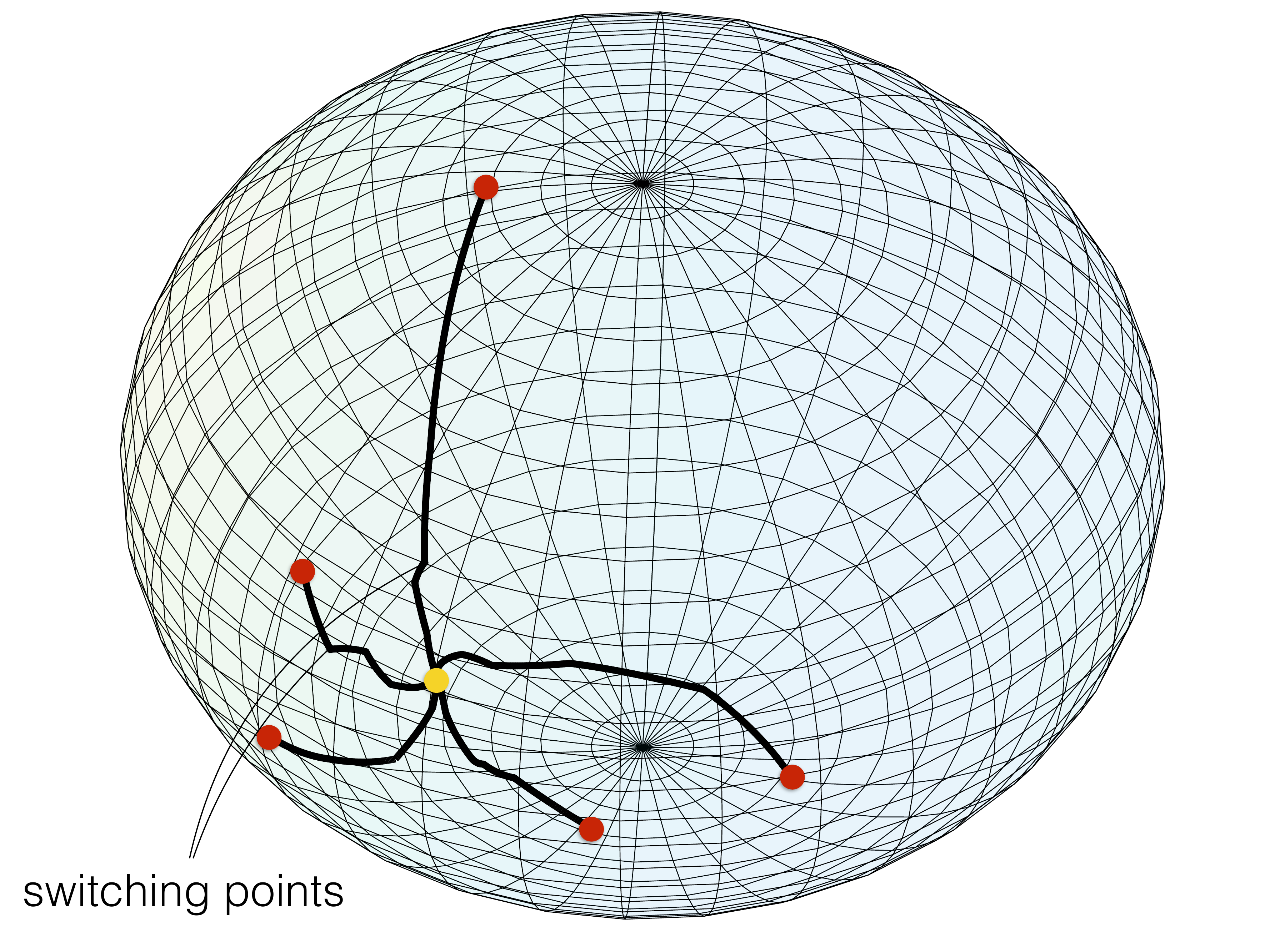}
  \caption{Convergence on a hemisphere.}
  \label{fig:10}
\end{figure}

To get a better understanding of Proposition~\ref{prop:10}, a numerical example is provided by Fig~\ref{fig:10}. In this example there are five agents with a uniformly quasi-strongly connected interaction graph. The agents were initially uniformly distributed on a hemisphere and the $f_{ij}$-functions were chosen to be constant; either equal to $1$ or $2$. In the figure, the red discs denote the initial positions and the yellow disc denote the final consensus point. We have also denoted two points on the trajectories where the graph switches. 

Now, the case when $0 \in \mathbb{R}^d$ is contained in the convex hull of the $z_i$'s is more intriguing. We provide the following result.

\begin{prop}\label{prop:4}
Suppose the dynamics for the $z_i$'s is governed by \eqref{eqq:11} and $\mathcal{G}(t)$ is uniformly strongly connected. Suppose $0 \in \mathbb{R}^d$ is contained in the convex hull of the $z_{i0}$'s, i.e., in $H(0,z_0)$,  and there is no $i$ such that $z_{i0} = 0$. Furthermore, suppose that the $z_{i0}$'s are contained in a compact set where the $f_{ij}$'s are 
bounded from below by a positive constant $K_d$.
Then the  set
$\mathcal{A}_z$---defined as the set where all the $z_i$'s are equal in the convex hull of the $z_{i0}$'s---is attractive. Furthermore, there is a fixed point $\bar z \in \mathbb{R}^d$ that all the $z_i$'s converge to.
\end{prop}

\emph{Proof:} \quad
We need to prove that the $z_i$'s converge to $\mathcal{A}_z$. 
In light of Proposition~\ref{prop:7}, the only case left to consider is when the point $0 \in \mathbb{R}^d$ is contained in the convex hull of the $z_i(t)$'s for all $t >0$, i.e., it is contained in $H(t,z_0)$ for all $t > 0$. We will thus only consider this case in the following where we need to prove that all the $z_i$'s converge to $0$. We partition this case into two sub-cases: \\ \\
\textbf{1)} the omega-limit set, denoted by $\Omega(z_0)$, does not contain a point $\bar z = [\bar z_1^T, \bar z_2^T, \ldots, \bar z_n^T]^T$ for which a $\bar z_i = 0$. \\
\noindent
\textbf{2)} the omega-limit set $\Omega(z_0)$ contains at least one point $\bar z = [\bar z_1^T, \bar z_2^T, \ldots, \bar z_n^T]^T$ where at least one of the $\bar z_i$'s is equal to zero. \\

\textbf{We begin by considering 1)}. There must be a ball $B$ around the origin such that there is no time $t$ for which a $z_i(t)$ is contained in the ball. This is proven in the following way. Proposition~\ref{prop:3} guarantees that no $z_i(t)$ can reach the origin in finite time. Thus, at any finite time $t_f$ there exists a largest open ball with radius $\epsilon(t_f)$ such that no $z_i(t)$ is contained in the ball during $[0,t_f]$. Assume that $\lim_{t_f \rightarrow \infty}\epsilon(t_f) = 0$. This implies that there is a point $\bar z = [\bar z_1^T, \bar z_2^T, \ldots, \bar z_n^T]^T$ that is in the closure of $\Omega(z_0)$ for which one of the $\bar z_i$'s is equal to zero. But $\Omega(z_0)$ is compact, hence we know that such a point also will be contained in $\Omega(z_0)$. This is a contradiction to the statement that such points are not contained in $\Omega(z_0)$. 

Now, in the set $B^n$ we replace the weights  
$f_{ij}(z_i,z_j)\frac{\|z_i\|}{\|z_j\|}$ in the right-hand side of \eqref{eqq:11}
by functions $h_{ij}(z_i,z_j)$ such that the total weights,
consisting of $f_{ij}(z_i,z_j)\frac{\|z_i\|}{\|z_j\|}$ outside $B^n$ and $h_{ij}(z_i,z_j)$ inside $B^n$, are 
globally Lipschitz on a set containing $(H(0,z_0))^n$
in the interior. Furthermore, $h_{ij}(y_i,y_j)$ is chosen to be positive when $y_i \neq y_j$. 

Now let us study the solution starting at $z_0$ at time $0$ of this modified system with the replaced weights in $B^n$. At any finite time the solution is the same as the original system. However, we can use the results in \cite{Johan_lyap_2017} to show that the solution to the modified system converges to $\mathcal{A}_z$ and in particular all the $z_i$'s converge to a fixed point that is nonzero. This means that after some
finite time $T$, the point $0$ is not contained in $H(t,z_0)$ for the modified as well as for the original system. But this is a contradiction to our assumption that the point $0 \in \mathbb{R}^d$ is contained in the convex hull of the $z_i(t)$'s for all $t >0$.\\ \\

\noindent 
\textbf{Now we consider 2)}. Let us first introduce $L(z) = \max_{i \in \mathcal{V}}\|z_i\|$. $L(z(t))$ is, besides continuous, monotonically decreasing. We assume that $\lim_{t \rightarrow \infty}L(z(t)) =\bar L> 0$. This means that for any $t \geq 0$, it 
holds that the set $\{j: \|z_j(t)\| \geq \bar L\}$ is nonempty. 
We will show that this assumption leads to a contradiction in the 
end of the proof. This, in turn, means that all the $z_i$'s converge to $0$.

Since the continuous $f_{ij}$'s are defined on a compact set, there is $K_u >0$ such that the $f_{ij}$'s are bounded from above by $K_u$. 
Furthermore, there is also an $M > 0$ such that $\|z_i(t)\| \leq M$ for all $i$ and $t \geq 0$. Take for example $M = L(z_0)$.

We continue by formulating a series of claims, each of which is followed by a proof. After these claims have been introduced, they are used as building blocks in the final part of the proof. Roughly, the claims can be understood as follows. The first claim says that if a $z_i$ is close to the origin, it will remain so for a specified time interval. The second claim says that if a $z_i$ has a neighbor that is close to the origin, it will will be ``dragged'' to the origin by this neighbor. The third claim simply says that there must be a $z_i$ close to the origin at some time. Then we show that that the $z_i$ that is close to the origin will drag all the other states close to the origin; so close that their distances to the origin is smaller than $\bar L$, which, in turn, is a contradiction.   \\

\noindent \textbf{Claim 1:} There is $\epsilon > 0$ satisfying $\epsilon^{\frac{1}{2}} \ll 1$ and $\epsilon^{\frac{1}{2}} \ll \bar L$ such that for time $\bar t \geq 0$ 
there is $i$ such that $z_i(\bar t)$ has smaller norm than $\epsilon$, then $\|z_{{i}}(t)\| \leq \epsilon^{\frac{1}{2}}$ for $t \in [\bar t, \bar t + (n+1)(T + \tau_D)]$, where $\tau_D$ is the lower bound on the dwell time and $T$ is the length of the time interval such that the union graph $\mathcal{G}([t, t + T))$ is guaranteed to be strongly connected, see Section~\ref{sec:preliminaries}. \\

It is assumed that $\epsilon^{\frac{1}{2}} \ll 1$ and $\epsilon^{\frac{1}{2}} \ll \bar L$. Suppose there is $\bar i$ and $\bar t$ such that $\|z_{\bar i}(\bar t)\| \leq \epsilon$.  
Let us consider the dynamics for $\|z_{\bar i}\|^2$. It is
\begin{align*}
 \frac{d}{dt}\|z_{\bar{i}}\|^2 & = 2\sum_{j \in \mathcal{N}_{\bar i}(t)}f_{ij}\frac{\|z_{\bar i}\|}{\|z_{\smash{j}}\|}(z_{\bar i}^Tz_j - \|z_{\bar i}\|^2) \\
& \leq 2\sum_{j \in \mathcal{N}_{\bar i}(t)}f_{ij}\frac{\|z_{\bar i}\|}{\|z_{\smash{j}}\|}z_{\bar i}^Tz_j \leq 2\sum_{j \in \mathcal{N}_{\bar i}(t)}K_u\|z_{\bar i}\|^2  \\
& \leq 2nK_u\|z_{\bar i}(t)\|^2  
\end{align*} 
Now we can use the Comparison Lemma to deduce that $\|z_{\bar i}(t)\| \leq \epsilon e^{nK_u(t -\bar t)}$ for $t > \bar t$. Now, if $\epsilon$ is  sufficiently small, the expression $\epsilon e^{nK_u(t -\bar t)}$ will be smaller than $\epsilon^{\frac{1}{2}}$ during $[\bar{t}, \bar{t} + (n+1)(T + \tau_D)]$. \\

\noindent \textbf{Claim 2:} There is $\epsilon > 0$ satisfying $\epsilon^{\frac{1}{3}} \ll 1$, $\epsilon^{\frac{1}{3}} \ll \bar L$ such that for any time $\bar t \geq 0$, if an agent $\bar i$ has a neighbor $\bar j$ during $[\bar t, \bar t + T]$ for which it holds that $\|z_{\bar j}(\bar t)\| \leq  \epsilon$, then there is a time $\bar t_2$ during $[\bar t, \bar t + T]$ such that $\min_{t \in [\bar t_2, \bar t_2 + \tau_D]}\|z_{\bar i}(t)\| \leq \max\{\epsilon^{\frac{1}{3}},Me^{{-\alpha \tau_D}}\}$  where $\alpha(\epsilon) = (K_d + nK_u) -  K_u\frac{\epsilon^{\frac{1}{3}}}{{\epsilon^{\frac{1}{2}}}}$. \\

Suppose that there is no $\bar t_2 \in [\bar t, \bar t + T]$ such that $\min_{t \in [\bar t_2, \bar t_2 + \tau_D]}\|z_{\bar i}(t)\| \leq \epsilon^{\frac{1}{3}}$. 
There is $\bar t_3 \in [\bar{t}, \bar t + T]$ such that $\bar{j} \in \mathcal{N}_{\bar i}(t)$ during the time interval $[\bar{t}_3, \bar t_3 +\tau_D] \subset [\bar{t}, \bar t + T + \tau_D]$, see the assumptions on the graph $\mathcal{G}(t)$ in Section~\ref{sec:preliminaries}. We assume that $\epsilon$ is small enough so that $\|z_{\bar{j}}(t)\| \leq \epsilon^{\frac{1}{2}}$ for $t \in [\bar t, \bar t + (n+1)(T + \tau_D)]$, see Claim 1. 

During the time interval $[\bar t_3, \bar t_3 + \tau_D]$ it holds that 
\begin{align}
\nonumber
& \frac{d}{dt}\|z_{\bar{i}}\|^2 = 2\sum_{j \in \mathcal{N}_{\bar i}(t)}f_{\bar ij}\frac{\|z_{\bar i}\|}{\|z_{\smash{j}}\|}(z_{\bar i}^Tz_j - \|z_{\bar i}\|^2) \\
\nonumber
& \leq 2f_{\bar i \bar j}\|z_{\bar i}\|^2(1 - \frac{\|z_{\bar i}\|}{\epsilon^{\frac{1}{2}}}) \\
\label{eq:new:3}
& \quad +2\sum_{j \in \mathcal{N}_{\bar i}(t)\backslash \bar j}f_{\bar ij}\frac{\|z_{\bar i}\|}{\|z_{\smash{j}}\|}(z_{\bar i}^Tz_j - \|z_{\bar i}\|^2).
\end{align}
Now we take a look at the first expression on the right-hand side of ``$\leq$'' in \eqref{eq:new:3}. We use the fact that $\|z_{\bar i}\| \geq \epsilon^{\frac{1}{3}}$ to obtain 
\begin{align*}
2f_{\bar i \bar j}\|z_{\bar i}\|^2(1 - \frac{\|z_{\bar i}\|}{\epsilon^{\frac{1}{2}}}) & \leq 2f_{\bar i \bar j}\|z_{\bar i}\|^2(1 - \frac{\epsilon^{\frac{1}{3}}}{{\epsilon^{\frac{1}{2}}}}) \\
& \leq 2K_d\|z_{\bar i}\|^2(1 - \frac{\epsilon^{\frac{1}{3}}}{{\epsilon^{\frac{1}{2}}}}).
\end{align*} 
Now we take a look at the (last) sum expression on the right-hand side of ``$\leq$'' in \eqref{eq:new:3}. For any $j$ such that $\|z_j\| \leq \epsilon^{\frac{1}{3}}$ the corresponding expression 
\begin{equation}\label{eq:new:2}
f_{\bar ij}\frac{\|z_{\bar i}\|}{\|z_{\smash{j}}\|}(z_{\bar i}^Tz_j - \|z_{\bar i}\|^2)
\end{equation}
is negative. For any $j \in \mathcal{N}_i(t)$ such that the expression in \eqref{eq:new:2} is positive, the following must hold
 \begin{align*}
& f_{\bar ij}\frac{\|z_{\bar i}\|}{\|z_{\smash{j}}\|}(z_{\bar i}^Tz_j - \|z_{\bar i}\|^2) \leq K_u\|z_{\bar i}\|^2
\end{align*} 

By using the above inequalities, we can conclude that 
\begin{align*}
& \frac{d}{dt}\|z_{\bar{i}}\|^2 \leq 2(K_d + nK_u)\|z_i(t)\|^2 - 2K_u\frac{\epsilon^{\frac{1}{3}}}{{\epsilon^{\frac{1}{2}}}}\|z_i(t)\|^2
\end{align*}
during the time interval $[t_3, t_3 + \tau_D]$. 
By using the Comparison Lemma we can conclude that 
\begin{align*}
& \|z_{\bar i}(\bar t_2 + \tau_D)\|^2 \leq \|z_{\bar i}(\bar t_2)\|^2 e^{-2\alpha \tau_D} \leq  M^2e^{-2\alpha \tau_D},
\end{align*}
where $\alpha = (K_d + nK_u) -  K_u\frac{\epsilon^{\frac{1}{3}}}{{\epsilon^{\frac{1}{2}}}}$. \\

\noindent \textbf{Claim 3:} For any $\epsilon > 0$ there is $\bar i \in \mathcal{V}$ and a corresponding
$\bar t$ such that $\|z_{\bar{i}}(\bar t)\| \leq \epsilon$. \\

Claim 3 is a consequence of the fact that the Omega limit set $\Omega(z_0)$ contains at least one point $\bar z = [\bar z_1^T, \bar z_2^T, \ldots, \bar z_n^T]^T$ where at least one of the $\bar z_i$'s is equal to zero. Thus there is $\bar i \in \mathcal{V}$ and a corresponding unbounded sequence $\{\bar t_n\}$ such that $\lim_{n \rightarrow \infty }\|z_{\bar i}(\bar t_n)\| = 0$.\\

Now, we can use the three claims above to obtain the sought contradiction concerning $L$'s convergence to $\bar L > 0$. Due to Claim 3, for any $\epsilon_1 > 0$ we know that there is a time $\bar t$ and an $\bar i_1$ for which $\|z_{\bar i_1}(\bar t)\| \leq \epsilon_1$. 

Let $\epsilon_2(\epsilon_1) = \max\{\epsilon_1^{\frac{1}{3}},Me^{{-\alpha(\epsilon_1) \tau_D}}\}$, where $\alpha$ is defined in Claim 2. By choosing $\epsilon_1$ small enough we can ensure, due to Claim 1, that $\epsilon_2^{\frac{1}{2}} \ll 1$ and $\epsilon_2^{\frac{1}{2}} \ll \bar L$ and there is $\bar i_2$ such that for $\bar i_1$, and $\bar i_2$, the norms of both $z_{\bar i_1}$ and $z_{\bar i_2}$ are smaller than $\epsilon_2^{\frac{1}{2}}$ during $[\bar{t} + T + \tau_D, \bar{t} + (n+1)(T + \tau_D)]$. 

Now let $\epsilon_3(\epsilon_2) = \max\{\epsilon_2^{\frac{1}{3}},Me^{{-\alpha(\epsilon_2) \tau_D}}\}$. Since $\epsilon_2$ is a function of $\epsilon_1$, by choosing $\epsilon_1$ small enough we can ensure, due to Claim 1, that $\epsilon_3^{\frac{1}{2}} \ll 1$ and $\epsilon_3^{\frac{1}{2}} \ll \bar L$ and there is $\bar i_3$ such that for $\bar i_1$, $\bar i_2$, and $\bar i_3$,the norms of $z_{\bar i}$, $z_{\bar i_2}$, and $z_{\bar i_3}$ are smaller than $\epsilon_3^{\frac{1}{2}}$ during $[\bar{t} + 2(T + \tau_D), \bar{t} + (n+1)(T + \tau_D)]$. 

By continuing in this manner, one can finally show that there is an $\epsilon_n < \bar L$ such that for all the $i$'s it holds that $\|z_i\| < \bar L$ at time $\bar{t} + (n+1)(T + \tau_D)$. But $L(t)$ is monotonically decreasing to $\bar L$ from above. Hence we have two contradictory statements.
\hfill $\blacksquare$

The point $0 \in \mathbb{R}^d$, plays a crucial role for this lifting method. If the $z_i$'s converge to a point that is not equal to $0$, then the $y_i$'s converge to a consensus formation. On the other hand, we do not know when convergence to $0$ for the $z_i$'s imply non-convergence to a consensus formation for the $y_i$'s.  

It has recently been shown that if the graph $\mathcal{G}(t)$ is static (or time-invariant) and symmetric, the $f_{ij}$'s fulfill certain differentiability assumptions, and $f_{ij} = f_{ji}$ for all $i,j$, then $\mathcal{A}$ is almost globally attractive under \eqref{eqq:6}~\citep{markdahl2017almost} (a simple choice of such $f_{ij}$'s  is when $f_{ij} = f_{ji} = a_{ij} = a_{ji} > 0$, i.e., the $a_{ij}$'s are positive scalars). On the other hand, the result does not hold for dimension $d = 2$.    
This means that under the same conditions on the graph and the $f_{ij}$'s as in \cite{markdahl2017almost}, the region of attraction of the point $0$ has measure zero when $d \geq 3$ and has positive measure when $d = 2$. Linking those results to the lifting method, e.g., by means of a geometric interpretation, remains an open problem.

\section{Conclusions}
\noindent
This paper addresses distributed synchronization or consensus on the unit sphere. A large 
class of consensus control laws is considered in which only relative information is used between neighbors for a time-switching interaction graph. We investigate how a lifting method can be used in the convergence analysis for these control laws. 
The proposed method is new in this context. It lifts the states from $\mathbb{S}^{d-1}$ to $\mathbb{R}^{d}$. In the higher-dimensional space, the dynamics of the states is described by a consensus protocol, where each agent is moving in a conical combination of the directions to its neighbors. The weights in the conical combination contain rational functions of the norms of the agents' states. 

The paper provides more general convergence results than has been reported before for hemispheres and furthermore provides convergence results for a consensus protocol with rational weights of the norms of the states. However, an additional purpose of the paper was---by introducing the lifting method---to hopefully serve as a stepping-stone towards future research on global convergence results for the considered consensus control laws.

\bibliographystyle{agsm}        
\bibliography{library.bib,ref.bib}

\end{document}